# On Pisier's construction of a polynomially bounded operator not similar to a contraction

John E. M$^c$Carthy *

M.S.R.I. and Washington University

March 12 1996

**Abstract**

A classical proof of Pisier's construction of a polynomially bounded operator not similar to a contraction is given.

## 0  Introduction

In [7], G. Pisier gave a very ingenious example of an operator that is polynomially bounded but not similar to a contraction, answering a question posed by Halmos in [4]. His proof uses martingales. The object of this expository note is to describe Pisier's example without using techniques of probability theory. We remark that this has also been done by S. Kisliakov and by K. Davidson and V. Paulsen.

Let $\mathcal{H}$ be a separable, infinite dimensional Hilbert space. Let $H^2(\mathcal{H})$ denote the space of $\mathcal{H}$-valued analytic functions $F$ on the disk for which

$$\|F\|^2 = \lim_{r \to 1} \frac{1}{2\pi} \int \|F(re^{i\theta})\|^2_{\mathcal{H}} d\theta = \frac{1}{2\pi} \int \|F(e^{i\theta})\|^2_{\mathcal{H}} d\theta < \infty$$

Let $K^2(\mathcal{H})$ denote the corresponding space of conjugate-analytic functions (not required to vanish at the origin). By passing to boundary values on the circle $\mathbb{T}$, both of these spaces can be thought of as subspaces of $L^2(\mathcal{H})$, the square-summable functions from $\mathbb{T}$ to $\mathcal{H}$. Let $\mathbf{Q}$ denote the orthogonal projection from $L^2(\mathcal{H})$ onto $K^2(\mathcal{H})$.

*The author was partially supported by the National Science Foundation grant DMS 9296099. Work done at M.S.R.I. is supported in part by NSF grant DMS-9022140.

Let $\{C_k\}_{k=1}^{\infty}$ be a fixed set of operators on $\mathcal{H}$ satisfying the canonical anti-commutation relations, *i.e.*

$$\begin{aligned} C_i C_j + C_j C_i &= 0 \\ C_i^* C_j + C_j C_i^* &= \delta_{ij} \end{aligned} \tag{0.1}$$

It is easy to construct $n$ such operators on a $2^n$ dimensional Hilbert space. For example, let

$$D = \begin{pmatrix} 0 & 1 \\ 0 & 0 \end{pmatrix} \quad E = \begin{pmatrix} 1 & 0 \\ 0 & -1 \end{pmatrix}$$

On the $n$-fold (or infinite) tensor product of $\mathbb{C}^2$ with itself, one can define

$$\begin{aligned} C_1 &= D \otimes I \otimes I \otimes \ldots \otimes I \\ C_2 &= E \otimes D \otimes I \otimes \ldots \otimes I \\ C_3 &= E \otimes E \otimes D \otimes \ldots \otimes I \\ &\ldots \\ C_n &= E \otimes E \otimes E \otimes \ldots \otimes E \otimes D \end{aligned}$$

Let $\Phi(z)$ be the $B(\mathcal{H})$-valued conjugate analytic function

$$\Phi(z) = \sum_{k=1}^{\infty} 2^{-k} \bar{z}^{2^k - 1} C_k$$

and let $\Gamma_\Phi : H^2(\mathcal{H}) \to K^2(\mathcal{H})$ be the vectorial Hankel operator

$$\Gamma_\Phi F = \mathbf{Q}(\Phi(F)).$$

Let $S : H^2(\mathcal{H}) \to H^2(\mathcal{H})$ be the operator of multiplication by $z$, and $X : K^2(\mathcal{H}) \to K^2(\mathcal{H})$ be multiplication by $z$ followed by the projection $\mathbf{Q}$. Finally, let $T$ be the operator on $K^2(\mathcal{H}) \oplus H^2(\mathcal{H})$ given by

$$T = \begin{pmatrix} X & \Gamma_\Phi \\ 0 & S \end{pmatrix} \tag{0.1}$$

Pisier proved that $T$ is polynomially bounded but not similar to a contraction.

Notice that, if $p$ is a polynomial and $M_p$ denotes multiplication by $p$, then

$$p(T) = \begin{pmatrix} p(X) & \Gamma_\Phi M_{p'} \\ 0 & p(S) \end{pmatrix} \tag{0.2}$$

(this can be checked for $p(z) = z^n$, and then follows for general $p$ by linearity).

# 1   $T$ is not similar to a contraction

First we need a lemma. The proof here is from [8, p. 65].

**Lemma 1.1** *Let $\alpha_1, \alpha_2, \ldots$ be complex numbers. Then*

$$\|\sum \alpha_i C_i\|^2 = \sum |\alpha_i|^2$$

PROOF: Let $R = \sum \alpha_i C_i$. Notice that, by (0.1),

$$\begin{aligned} R^*R + RR^* &= (\sum \bar{\alpha}_i C_i^*)(\sum \alpha_j C_j) + (\sum \alpha_j C_j)(\sum \bar{\alpha}_i C_i^*) \\ &= \sum_{i,j} \bar{\alpha}_i \alpha_j (C_i^* C_j + C_j C_i^*) \\ &= \sum |\alpha_i|^2 \end{aligned}$$

and

$$\begin{aligned} R^2 &= \frac{1}{2}(\sum_{i,j} \alpha_i \alpha_j C_i C_j + \sum_{i,j} \alpha_i \alpha_j C_j C_i) \\ &= 0 \end{aligned}$$

Therefore

$$\begin{aligned} \|R\|^4 &= \|R^* R R^* R\| \\ &= \|R^*(R^*R + RR^*)R\| \\ &= \|R^*(\sum |\alpha_i|^2)R\| \\ &= \sum |\alpha_i|^2 \|R\|^2 \end{aligned}$$

So $\|R\|^2 = \sum |\alpha_i|^2$, as required. $\square$

An operator $R$ is called completely polynomially bounded if there exists a constant $c$ such that for every matrix valued (or operator valued) polynomial $P$, the inequality $\|P(R)\| \leq c\|P\|$ holds, where $\|P\|$ is the supremum of the operator norms of $P(z)$ as $z$ ranges over the unit disk. If $R$ were similar to a contraction, then it would be completely polynomially bounded. Reason: suppose $R = A^{-1}BA$, and $B$ is a contraction. Let $P(z) = \sum z^n D_n$ be an operator valued polynomial. Then

$$P(R) = \sum R^n \otimes D_n = (A^{-1} \otimes I)P(B)(A \otimes I)$$

By the operator-valued version of von Neumann's inequality (which can be proved using the Sz.-Nagy dilation theorem just as the scalar one), $\|P(B)\| \leq \|P\|$, so $\|P(R)\| \leq$ $\|A^{-1}\|\|A\|\|P\|$. The converse of this statement - completely polynomially bounded implies similar to a contraction - is true, but we shall not use this fact [6].

We shall prove that $T$ is not completely polynomially bounded, and therefore cannot be similar to a contraction. Fix a positive integer $n$. Let

$$P(z) = \sum_{k=1}^{n} z^{2^k} C_k$$

By Lemma (1.1), we have $\|P\| = \sqrt{n}$. Let us calculate the norm of

$$P(T) = \sum_{k=1}^{n} T^{2^k} \otimes C_k$$

With respect to the 2-by-2 decomposition of $T$ in (0.1), this also has a 2-by-2 decomposition, with the $(1,2)$ entry being

$$R = \Gamma_\Phi \sum_{k=1}^{n} 2^k M_{z^{2^k}-1} \otimes C_k = \mathbf{Q} \sum_{k=1}^{n} 2^k (\sum_{j=1}^{\infty} 2^{-j} \bar{z}^{2^j-1} z^{2^k-1} C_j) \otimes C_k.$$

The constant term in $R$ is then $\sum_{k=1}^{n} C_k \otimes C_k$, and it will suffice to prove that the norm of this is larger than $c\sqrt{n}$ for any fixed $c$.

Let $\{\mathbf{e}_i\}_{i=1}^{\infty}$ be an orthonormal basis for $\mathcal{H}$, and assume that the $C_k$ are constructed as in the previous section, so in particular $\langle C_k \mathbf{e}_i, \mathbf{e}_j \rangle$ is always real (this assumption is to avoid introducing an antilinear isometry of $\mathcal{H}$). Then

$$\begin{aligned}
\|\sum_{k=1}^{n} C_k \otimes C_k\| &\geq \frac{1}{2^n} \langle \sum_{k=1}^{n} C_k \otimes C_k (\sum_{i=1}^{2^n} \mathbf{e}_i \otimes \mathbf{e}_i), (\sum_{j=1}^{2^n} \mathbf{e}_j \otimes \mathbf{e}_j) \rangle \\
&= \frac{1}{2^n} \sum_{i,j,k} \langle C_k \mathbf{e}_i, \mathbf{e}_j \rangle \langle C_k \mathbf{e}_i, \mathbf{e}_j \rangle = \frac{1}{2^n} \sum_{i,j,k} \langle C_k \mathbf{e}_i, \mathbf{e}_j \rangle \langle \mathbf{e}_j, C_k \mathbf{e}_i \rangle \\
&= \frac{1}{2^n} \sum_{k=1}^{n} \sum_{i=1}^{2^n} \langle C_k \mathbf{e}_i, C_k \mathbf{e}_i \rangle \\
&= \frac{1}{2^n} \sum_{k=1}^{n} 2^{n-1} = \frac{n}{2}
\end{aligned}$$

Therefore $\|p(T)\|/\|p\| \geq \sqrt{n}/2$ and is not bounded, so $T$ is not completely polynomially bounded, and hence is not similar to a contraction.

## 2 $T$ is polynomially bounded

To prove that $T$ is polynomially bounded, by (0.2) we must prove that there exists a constant $c$ such that $\|\Gamma_\Phi M_{p'}\| \leq c\|p\|_\infty$ for all polynomials $p$. Now,

$$\|\Gamma_\Phi M_{p'}\| = \sup\{|\langle \Gamma_\Phi p' F, G \rangle| : \|F\|_{H^2(\mathcal{H})} = 1, \|G\|_{K^2(\mathcal{H})} = 1\}$$

Let
$$F(z) = \sum_{n=0}^{\infty} z^n \xi_n, \quad G(z) = \sum_{m=0}^{\infty} \bar{z}^m \eta_m$$
where
$$\sum_{n=0}^{\infty} \|\xi_n\|_{\mathcal{H}}^2 = 1 = \sum_{m=0}^{\infty} \|\eta_m\|_{\mathcal{H}}^2.$$
Let $\sigma$ denote normalized Lebesgue measure on $\mathbb{T}$. Then
$$\langle \Gamma_\Phi p' F, G \rangle = \int_{\mathbb{T}} d\sigma(z) p'(z) \sum_{k=1}^{\infty} 2^{-k} \bar{z}^{2^k-1} \sum_{n,m=0}^{\infty} z^{n+m} \langle C_k \xi_n, \eta_m \rangle.$$

By (1.1) the map $v : \mathcal{H} \to B(\mathcal{H})$ which sends $\mathbf{e}_i$ to $C_i$ is an isometry. Therefore
$$\begin{aligned}
\langle C_k \xi_n, \eta_m \rangle &= tr[C_k(\xi_n \otimes \eta_m)] \\
&= (v(\mathbf{e}_k), \xi_n \otimes \eta_m)_{tr} \\
&= (\mathbf{e}_k, v^t(\xi_n \otimes \eta_m))_{tr}
\end{aligned}$$

where $(.,.)_{tr}$ is the linear pairing between $B(H)$ and the trace class operators, given by $(A, B)_{tr} = tr(AB)$. Therefore
$$\langle \Gamma_\Phi p' F, G \rangle = \int_{\mathbb{T}} d\sigma(z) p'(z) \sum_{k=1}^{\infty} 2^{-k} \bar{z}^{2^k-1} (v(\mathbf{e}_k), [\sum_{n=0}^{\infty} z^n \xi_n] \otimes [\sum_{m=0}^{\infty} z^m \eta_m])_{tr} \qquad (2.1)$$

Now, $[\sum_{n=0}^{\infty} z^n \xi_n] \otimes [\sum_{m=0}^{\infty} z^m \eta_m]$ is in $H^1(\mathcal{H} \otimes \mathcal{H})$, and moreover

$$\begin{aligned}
\|[\sum_{n=0}^{\infty} z^n \xi_n] \otimes [\sum_{m=0}^{\infty} z^m \eta_m]\|_{H^1(\mathcal{H} \otimes \mathcal{H})} &= \int_{\mathbb{T}} d\sigma(z) \|[\sum_{n=0}^{\infty} z^n \xi_n] \otimes [\sum_{m=0}^{\infty} z^m \eta_m]\|_{\mathcal{H} \otimes \mathcal{H}} \\
&= \int_{\mathbb{T}} d\sigma(z) \|\sum_{n=0}^{\infty} z^n \xi_n\|_{\mathcal{H}} \|\sum_{m=0}^{\infty} z^m \eta_m\|_{\mathcal{H}} \\
&\leq [\int_{\mathbb{T}} d\sigma(z) \|\sum_{n=0}^{\infty} z^n \xi_n\|_{\mathcal{H}}^2]^{1/2} [\int_{\mathbb{T}} d\sigma(z) \|\sum_{m=0}^{\infty} z^m \eta_m\|_{\mathcal{H}}^2]^{1/2} \\
&= \|F\|_{H^2(\mathcal{H})} \|G\|_{K^2(\mathcal{H})} \\
&= 1
\end{aligned}$$

Combining this with (2.1), we get that

$$\begin{aligned}
\|\Gamma_\Phi M_{p'}\| &\leq \|\mathbf{Q} \sum_{k=1}^{\infty} 2^{-k} p'(z) \bar{z}^{2^k-1} \mathbf{e}_k\|_{H^1(\mathcal{H})^*} \| \sum_{n,m=0}^{\infty} z^{n+m} v^t(\xi_n \otimes \eta_m)\|_{H^1(\mathcal{H})} \\
&\leq \|\mathbf{Q} \sum_{k=1}^{\infty} 2^{-k} p'(z) \bar{z}^{2^k-1} \mathbf{e}_k\|_{H^1(\mathcal{H})^*}
\end{aligned}$$

The dual of $H^1(\mathcal{H})$ is $BMO(\mathcal{H})$ (in [1], Bourgain proves that for any Banach space $X$, the dual of $H^1(X)$ is $BMO(X^*)$; however if $X$ is a Hilbert space, the proof of the duality follows from the same argument as is used in the scalar case - *i.e.* Fefferman's theorem - in *e.g.* [3]). Therefore we are done if we can prove the following:

$$\|\mathbf{Q}\sum_{k=1}^{\infty} 2^{-k} p'(z)\bar{z}^{2^k-1}\mathbf{e}_k\|_{BMO(\mathcal{H})} \leq c\|p\|_{\infty} \tag{2.2}$$

We do this in section 3.

## 3  The key estimate

Let

$$K_k(e^{i\theta}) := 2^{-k} \sum_{n=-2^k}^{2^k} (2^k - |n|)e^{in\theta}$$

be the Fejér kernel, and let

$$T_k f(e^{i\theta}) = f * K_k(e^{i\theta}) = \frac{1}{2\pi}\int_{\mathbb{T}} f(e^{i\phi})K_k(e^{i(\theta-\phi)})d\phi$$

Then what we want to estimate in (2.2) can be written

$$\mathbf{Q}\sum_{k=1}^{\infty} 2^{-k}p'(z)\bar{z}^{2^k-1}\mathbf{e}_k = \mathbf{Q}\sum_{k=1}^{\infty} T_k p(z)\bar{z}^{2^k}\mathbf{e}_k, \tag{3.1}$$

where we write $T_k p(z)\bar{z}^{2^k}$ to mean $T_k[p(z)\bar{z}^{2^k}]$, *i.e.*

$$[T_k p(z)\bar{z}^{2^k}](e^{i\theta}) = \frac{1}{2\pi}\int_{\mathbb{T}} p(e^{i\phi})e^{-i2^k\phi}K_k(e^{i(\theta-\phi)})d\phi$$

As was observed in this context by S. Kisliakov, the norm in the dual of $H^1(\mathcal{H})$ depends only on the conjugate analytic part of the function; so one can remove the $\mathbf{Q}$ in (3.1), and it is sufficient to show

$$\|\sum_{k=1}^{\infty} T_k p(z)\bar{z}^{2^k}\mathbf{e}_k\|_{BMO(\mathcal{H})} \leq c\|p\|_{\infty} \tag{3.2}$$

Now

$$\|\sum_{k=1}^{\infty} T_k p(z)\bar{z}^{2^k}\mathbf{e}_k\|_{BMO(\mathcal{H})}^2 = \left[\sup_I \frac{1}{|I|}\int_I \left[\sum_{k=1}^{\infty} |T_k p(z)\bar{z}^{2^k} - MV_k|^2\right]^{1/2} d\theta\right]^2$$

$$\leq \sup_I \sum_{k=1}^{\infty} \frac{1}{|I|}\int_I |T_k p(z)\bar{z}^{2^k} - MV_k|^2 d\theta \tag{3.3}$$

where $I$ ranges over the arcs of the circle, and $MV_k$ denotes the mean value, over $I$, of the preceding expression.

Let us fix some arc $I$; we shall estimate (3.3) by breaking the sum into two pieces, $k$ small and $k$ large. We shall let $c$ denote a universal constant, independent of $I, k$ and $p$, that may vary from one line to the next.

Case (a): $2^{-k} \geq |I|$.

We shall estimate the size of the derivative of $T_k p(z) \bar{z}^{2^k}$. Note that for any function $f$,

$$|\frac{d}{d\theta} T_k f(e^{i\theta})| \leq \|f\|_\infty \|\frac{d}{d\theta} K_k(e^{i\theta})\|_1$$

The Fejér kernel can be expressed in closed form as

$$K_k(e^{i\theta}) = \frac{1}{2^k} \left( \frac{\sin(2^{k-1}\theta)}{\sin(\theta/2)} \right)^2 \tag{3.3}$$

(see e.g. [5]), and one can estimate the $L^1$ norm of the derivative by summing the variation over intervals of length $\pi/2^k$, to get $c 2^k$. Therefore

$$\|\frac{d}{d\theta} T_k f\|_i \leq c 2^k \|f\|_\infty$$

Now let the end-points of $I$ be $e^{ia}$ and $e^{ib}$. Then for any function $f$

$$\begin{aligned}
\frac{1}{|I|} \int_I |T_k f - MV|^2 d\theta &\leq \frac{1}{|I|} \int_a^b |T_k f(e^{i\theta}) - T_k f(e^{ia})|^2 d\theta \\
&\leq \frac{1}{|I|} \int_a^b \|\frac{d}{d\theta} T_k f\|_\infty (\theta - a)^2 d\theta \\
&\leq c 2^{2k} |I|^2 \|f\|_\infty^2
\end{aligned}$$

Letting $f$ be $\bar{z}^{2^k} p(z)$, and summing over $k$ such that $2^k$ is less than or equal to $1/|I|$, we get

$$\sum_{2^{-k} \geq |I|} \frac{1}{|I|} \int_I |T_k p(z) \bar{z}^{2^k} - MV_k|^2 d\theta \leq c \|p\|_\infty^2 \tag{3.4}$$

Case (b): $2^{-k} \leq |I|$.

For large $k$, we can drop the mean value term. Let $3I$ denote the interval concentric with $I$ of three times the length. Write $p = p_1 + p_2$, where $p_1$ is supported on $3I$, and $p_2$ is supported on $\mathbb{T} \setminus 3I$.

Case (i): $p_1$.

Note that

$$\sum_{k=1}^\infty \int_\mathbb{T} |T_k \bar{z}^{2^k} p_1|^2 = \sum_{k=1}^\infty \{ \sum_{n=0}^{2^k} |\hat{p}_1(n)|^2 \frac{n^2}{2^{2k}} + \sum_{n=2^k+1}^{2^{k+1}} |\hat{p}_1(n)|^2 \frac{(2^{k+1} - n)^2}{2^{2k}} \}$$

Interchanging the order of summation, one gets that the sum is less than or equal to

$$\frac{7}{3}\sum_{n=0}^{\infty}|\hat{p}_1(n)|^2 \leq \frac{7}{3}\|p_1\|_2^2 \leq 7|I|\,\|p\|_\infty^2$$

Therefore
$$\sum_{k=1}^{\infty}\frac{1}{|I|}\int_I |T_k p_1(z)\bar{z}^{2^k} - MV_k|^2 d\theta \leq c\|p\|_\infty^2 \qquad (3.5)$$

Case (ii): $p_2$.

We have

$$\sum_{2^{-k}\leq |I|}\frac{1}{|I|}\int_I |T_k \bar{z}^{2^k} p_2|^2 d\theta = \sum_{2^{-k}\leq |I|}\frac{1}{|I|}\int_I d\theta \left|\int_{\mathbb{T}\setminus I\!I} d\phi\, e^{-i2^k\phi} p(e^{i\phi}) K_k(e^{i(\theta-\phi)})\right|^2$$

$$\leq \sum_{2^{-k}\leq |I|}\frac{1}{|I|}\int_I d\theta \|p\|_\infty^2 \left|\int_{|I|/2}^{2\pi-|I|/2} K_k(e^{i\phi}) d\phi\right|^2$$

From (3.3) we get that
$$\int_{|I|/2}^{2\pi-|I|/2} K_k(e^{i\phi}) d\phi \leq c 2^{-k}\frac{1}{|I|}.$$

Therefore
$$\sum_{2^{-k}\leq |I|}\frac{1}{|I|}\int_I |T_k \bar{z}^{2^k} p_2|^2 d\theta \leq \sum_{2^{-k}\leq |I|}\|p\|_\infty^2 2^{-2k}\frac{1}{|I|^2}$$
$$\leq c\|p\|_\infty^2 \qquad (3.6)$$

Combining (3.4), (3.5) and (3.6) we get

$$\sup_I \frac{1}{|I|}\sum_{k=1}^{\infty}\int_I |T_k p(z)\bar{z}^{2^k} - MV_k|^2 d\theta \leq c\|p\|_\infty^2$$

as desired.

# 4  Finite dimensional estimates

Let $R$ be any operator on an $n$-dimensional Hilbert space. There are two quantities associated with it (which may be infinite):

$$\|R\|_{s.b.} := \inf\{\|A\|\|A^{-1}\| \;:\; \|ARA^{-1}\| \leq 1\}$$
$$\|R\|_{p.b.} := \sup\{\|p(R)\| \;:\; p \text{ a polynomial and } \|p\|_\infty \leq 1\}$$

Adopt the convention that if $R$ is not similar to any contraction, then $\|R\|_{s.b.} = \infty$. By von Neumann's inequality, one always has $\|R\|_{p.b.} \leq \|R\|_{s.b.}$. In [2], J. Bourgain proved that, for any matrix $R$ on an $n$-dimensional Hilbert space,

$$\|R\|_{s.b.} \leq \|R\|_{p.b.}^4 \log n$$

One of the consequences of Pisier's result, as he points out, is that there is a lower bound:

**Corollary 4.1 (Pisier)** *There is a constant $\delta > 0$ such that, for all $K > 1$, and all $n \geq 2$, there exists an $n$-by-$n$ matrix $R$ with $\|R\|_{p.b.} \leq K$ and $\|R\|_{s.b.} \geq \delta(K-1)\sqrt{\log n}$.*

PROOF: It is sufficient to prove for $n \geq 8$ and a power of two; let $n = 2^{m+2}$. Let $\mathcal{H}_m = H^2(l_{2^m}^2) \ominus z^{2^m} H^2(l_{2^m}^2)$, and $\mathcal{K}_m = K^2(l_{2^m}^2) \ominus \bar{z}^{2^m} K^2(l_{2^m}^2)$. Let $C_1, \ldots, C_m$ satisfy the canonical anti-commutation relations on $l_{2^m}^2$, let $\Phi_m(z) = \sum_{k=1}^m 2^{-k} \bar{z}^{2^k-1} C_k$ and let $\Gamma_{\Phi_m} : \mathcal{H}_m \to \mathcal{K}_m$ be $\Phi_m$ followed by projection from $L^2(l_{2^m}^2)$ onto $\mathcal{K}_m$. Let $S_m$ and $X_m$ be truncated versions of the forward and backward shifts, and on $\mathcal{K}_m \oplus \mathcal{H}_m$, let

$$T_m = \begin{pmatrix} X_m & \varepsilon \Gamma_{\Phi_m} \\ 0 & S_m \end{pmatrix}$$

where $\varepsilon > 0$ will be chosen later. Note that $T$ is an $n$-by-$n$ matrix. Then for any polynomial $p$,

$$p(T_m) = \begin{pmatrix} p(X_m) & \varepsilon \Gamma_{\Phi_m} p'(S_m) \\ 0 & p(S_m) \end{pmatrix}$$

Just as above, one gets that

$$\|T_m\|_{s.b.} \geq \varepsilon \sqrt{m}/2$$

and $\|\Gamma_{\Phi_m} p'(S_m)\| \leq c\|p\|_\infty$ where $c$ is the constant in 2.2, and so

$$\|T_m\|_{p.b.} \leq \varepsilon c + 1.$$

Now let $\varepsilon = (K-1)/c$, and let $\delta = 1/5c$. Then $T_m$ works.

□

## Acknowledgements

The author would like to thank M. Christ and A. Arias for valuable discussions.